# Congestion Relief and Load Curtailment Reduction with FACTS Devices

Barati Masoud, Shayanfar Heidar Ali, Kazemi Ahad

Center of Excellent for Power System Automation Department, Iran University of science and Technology Tehran, Iran

*Abstract*—This paper investigates the effect of FACTS Devices on voltage stability constrained Optimal Power Flow (OPF) formulation where the objective function is to minimize power system load curtailment. Incorporating FACTS Devices affects the topology and hence the power flow distribution. Some numerical cases are presented to discuss the effects of FACTS Devices that is expected when it is incorporated into the load curtailment formulation. The test result reflects the impacts FACTS Devices has in reducing load curtailment during a line congestion situation. The paper also discusses the applicability of approach in security based reliability studies for systems having control components like FACTS Devices.

## I. INTRODUCTION

In deregulated power systems, economic competition leads to an understatement as regards to maintaining security features of overall system. One such security issue is the voltage stability of system. Several voltage instability incidents have been reported, in the recent past, all over the globe. These are results of operating system with very less voltage stability margin under normal conditions. Another offshoot of the deregulated operating environment is existence of various transaction paths based on the location of energy market players and their contract amount. Thus, congestion management has become one of an important operational issues. In a deregulated environment, congestion alleviation would mean load curtailment in certain situations. Authors have proposed an OPF based formulation, incorporating the voltage stability as an additional constraint, for evaluating the amount of load curtailment [1]. The voltage stability margin indicator used was discussed in the literature of an earlier reference [2]. The conventional OPF procedure has been well documented in many earlier works [3,4,5,6]. Further, the application of voltage stability constrained OPF procedure in evaluating security based composite system reliability was dealt with in a subsequent work [7]. Thus, in the emerging deregulation market any control action has to incorporate security features to maintain an acceptable level of system reliability.

It is observed that when the OPF results into a load Curtailment, some of the branches still have thermal capacity limits that are still under-utilized. Selective change of branch parameters can reduce the amount of load curtailment for a given loading scenario. Providing voltage support at weak voltage points can also help in reducing the load curtailments [8]. The above is the motivation for us in exploring the effects of Flexible AC Transmission System (FACTs) components into the voltage stability constrained OPF formulation.

FACTS Devices has been used to enhance angle stability and to mitigate the sub-synchronous resonance [9, 10, 11]. However, the use of TCSC in redistributing power flow and its effect on load curtailment has also been demonstrated in previous works [8, 12, and 13]. Static Var Compensator (SVC) used as reactive support at voltage weak points can also reduce potential load curtailment [8].

This paper first presents the method of including the steady state model for FACTS Devices into the voltage stability constrained OPF formulation. The impacts realized, by incorporating FACTS devices, are then brought about by case study. For comparison purposes the results as obtained without FACTS devices following the voltage stability constrained algorithm [1, 7] is also evaluated. Moreover, it also discusses the applicability of the above methodology in a security based composite reliability assessment.

## II. THE MODEL OF FACTS DEVICES FOR FLOW CONTROL

The FACTS model used in this article is the PIM [14, 15,16]. This model can keep the symmetric characteristic of the admittance matrix and also adapt to different types of FACTS devices. The configuration of the Unified Power Flow Controller (UPFC) and the development of its PIM model are summarized in Figure 1 and 2, respectively. In constraints, how to integrate the effects of the viewpoint of power markets and the control of FACTS devices with the OPF is proposed in this paper as follows. For any nodes $i$ without the FACTS devices between node $i$ and another node:

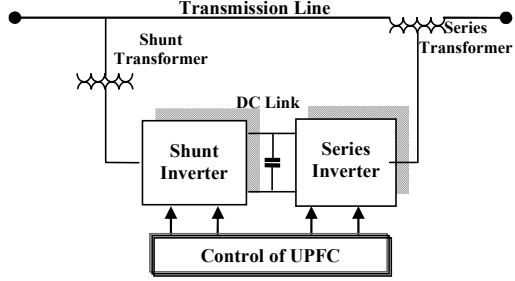

Fig. 1. Configuration of the UPFC

$$P_i - P(V,\theta) = 0 \quad (1)$$

$$Q_i - Q(V,\theta) = 0 \quad (2)$$

For any FACTS devices installed between node $i$ and node $j$:

$$P_I + \beta_1 P_{I(inj)} - P(V,\theta) = 0 \quad (3)$$

$$Q_I + \beta_1 Q_{IL(inj)} + \beta_2 Q_{II(inj)} - Q(V,\theta) = 0 \quad (4)$$

$$P_J + \beta_1 P_{J(inj)} - P(V,\theta) = 0$$

$$Q_J + \beta_1 Q_{J(inj)} - Q(V,\theta) = 0 \quad (5)$$

Where $\beta_1 = \beta_2 = 1$ for a UPFC, $\beta_1 = 0$ and $\beta_2 = 1$ for a

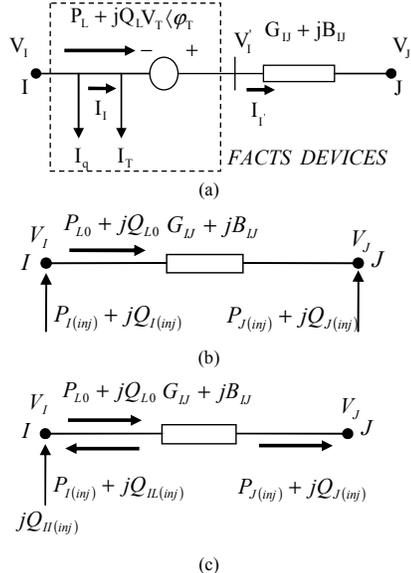

(c)
Figure 2. Development of the decomposed PIM: (a) voltage-current source model of the FACTS devices, (b) conventional PIM of the FACTS devices, (c) decomposed PIM of the FACTS devices.

Static Synchronous Compensation (STATCOM), $\beta_1 = 1$ and $\beta_2 = 0$ for a Static Synchronous Series Compensation (SSSC).

For any FACTS devises installed between node $I$ and node $J$ the following apply. Active power injection limit:

$$P_{I(inj)}^{min} \leq P_{I(inj)} \leq P_{I(inj)}^{max} \quad (6)$$

Reactive power injection limit:

$$Q_{IL(inj)}^{min} \leq Q_{IL(inj)} \leq Q_{IL(inj)}^{max} \quad (7)$$

$$Q_{II(inj)}^{min} \leq Q_{II(inj)} \leq Q_{II(inj)}^{max} \quad (8)$$

Where (19) For UPFC and SSSC, (20) for UPFC and (21) for UPFC and STATCOM.

To calculated $Q_{J(inj)}$:

$$\vec{V}_I' = \vec{V}_I + \vec{V}_T \quad (9)$$

$$Arg(\vec{I}_q) = Arg(\vec{V}_I) \pm 90^0 \quad (10)$$

$$\vec{I}_T = \frac{\text{Re}(\vec{V}_T \vec{I}_I^*)}{\vec{V}_I} \quad (11)$$

Thus,

$$\vec{S}_{I(inj)} = \vec{V}_I \left[ -\vec{I}_T - \vec{I}_q - \vec{V}_T (G_{IJ} + jB_{IJ}) \right]^* \quad (11)$$

$$P_{I(inj)} = G_{IJ} \begin{bmatrix} -V_T^2 - 2V_I V_T \cos(\theta_I - \varphi_T) \\ + V_J V_T \cos(\theta_J - \varphi_T) \end{bmatrix} \\ - B_{IJ} [V_J V_T \sin(\theta_J - \varphi_T)] \quad (12)$$

$$Q_{IL(inj)} = V_I V_T \begin{bmatrix} G_{IJ} \sin(\theta_I - \varphi_T) \\ + B_{IJ} \cos(\theta_I - \varphi_T) \end{bmatrix} \quad (13)$$

$$Q_{II(inj)} = -V_I V_q \quad (14)$$

$$\vec{S}_{J(inj)} = \vec{V}_J \left[ \vec{V}_T (G_{IJ} + jB_{IJ}) \right]^* \quad (15)$$

$$P_{J(inj)} = V_J V_T \begin{bmatrix} G_{IJ} \cos(\theta_J - \varphi_T) \\ + B_{IJ} \sin(\theta_J - \varphi_T) \end{bmatrix} \quad (16)$$

$$Q_{J(inj)} = V_J V_T \begin{bmatrix} G_{IJ} \cos(\theta_J - \varphi_T) \\ - B_{IJ} \sin(\theta_J - \varphi_T) \end{bmatrix} \quad (17)$$

The modified OPF can be solved by a nonlinear programming technique that takes into account
- Active and reactive power
- Load Curtailment Reduction
- FACTS devices

With the proposed method, the ISO would operate the transmission network more effectively. The power flow can be controlled by FACTS devices and coordinated with the transaction trade simultaneously.

## III. VOLTAGE STABILITY CONSTRAINED LOAD CURTAILMENT FORMULATION WITH FACTS DEVICES

The formulation for incorporating the FACTS devices control, into the procedure proposed by the authors [1], is presented herewith. The detail of the voltage stability index, L is described in published literatures [1, 8].

To simplify the simulations we have kept the load power factor to be constant i.e. we assume that when a certain amount of real load has been shed at one bus, the corresponding reactive load will also be shed.

It can be observed in the OPF formulation that it includes

Objective $min \sum_{i=1}^{n} \Delta P_{li}$

S.T:

Active power equation with FACTS devices at bus $i$ (18)
Reactive power equation with FACTS devices at bus $i$ (19)

$$P_{li}/P_{lireq} = Q_{li}/Q_{lireq} \quad (20)$$
$$0 \le P_{li} \le P_{lireq} \quad (21)$$
$$0 \le Q_{li} \le Q_{lireq} \quad (22)$$
$$|V_i|_{min} \le |V_i| \le |V_i|_{max} \quad (23)$$
$$P_{gi\,min} \le P_{gi} \le P_{gi\,max} \quad (24)$$
$$Q_{gi\,min} \le Q_{gi} \le Q_{gi\,max} \quad (25)$$
$$P_{ij}^2 + Q_{ij}^2 \le S_{ij}^2 \quad (26)$$
$$L_i \le L_{crit} \quad (27)$$
equations(3) to (17) (28)

Here,
$\Delta P_{li} = P_{lireq} - P_{li}$

Where,

$P_{lireq}$ : Real load demand at bus $i$

$P_{li}$ : Actual real load supply at bus $i$

$n$ : Total number of load flow buses in the system

$P_{gi}$ : Real power generation at bus $i$

$Q_{gi}$ : Reactive power generation at bus $i$

$Q_{lireq}$ : Reactive load demand at bus $i$

$Q_{li}$ : Actual reactive load supply at bus $i$

$|V_i|$ : Voltage magnitude at bus $i$

$G_{ij}, B_{ij}$ : Real/reactive part of the $ij^{th}$ element of bus admittance matrix

$L_i$ is the index $L$ evaluated at the $i^{th}$ bus other than the generation buses

$L_{crit}$ is the threshold value of the index acceptable for the system

the power balance equations (18, 19) generation limits (24, 25), line loading limits (26), voltage magnitude limits (23) and voltage stability constraint (27).

For the load curtailment policy, which we have adopted, i.e. constant power factor, an additional constraint (20) has been added. This works along with the allowable range of the possible load that can be supplied, which is represented by equalities (21) and (22). To incorporate the FACTS devices control into the OPF's description the constraint (28) has been introduced.

## IV. THE EFFECT OF FACTS DEVICES ON LOAD CURTAILMENT

### A. Case Study

We shall now apply the OPF to study the WSCC-9 bus test system as shown in Fig.3. The line parameters and the thermal limits are as shown in Table I.

Three studies were carried out for the above test system. The details of the focus, the simulation set-up conditions and

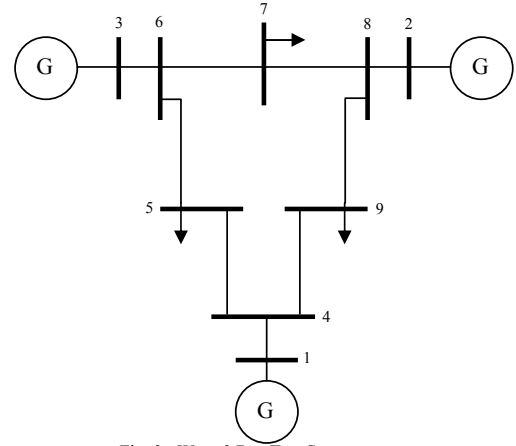

Fig. 3. Wscc 9 Bus Test System

TABLE I

LINE PARAMETERS AND LOADING LIMITS

| Line | Resistance (p.u) | Reactance (p.u) | Susceptance (p.u) | MVA Rating |
|------|------------------|-----------------|-------------------|------------|
| 1-4  | 0.0000           | 0.0576          | 0.0000            | 250        |
| 4-5  | 0.0170           | 0.0920          | 0.1580            | 250        |
| 5-6  | 0.0390           | 0.1700          | 0.3580            | 150        |
| 3-6  | 0.0000           | 0.0586          | 0.0000            | 300        |
| 6-7  | 0.0119           | 0.1008          | 0.2090            | 150        |
| 7-8  | 0.0085           | 0.0720          | 0.1490            | 250        |
| 8-2  | 0.0000           | 0.0625          | 0.0000            | 250        |
| 8-9  | 0.0320           | 0.1610          | 0.3060            | 250        |
| 9-4  | 0.0100           | 0.0850          | 0.1760            | 250        |

their details are discussed in the following sub-sections.

1) Before evaluating the effect of FACTS devices incorporation into the system, let us first demonstrate the choice of $L_{crit}$ value for the test system. Load bus 5 was supposedly having a load demand of 90+ j 30 MVA, bus 7 a load demand of 100 + j 30 MVA and load bus 9 having demand of 149.67 + j 59.87 MVA. All the generator buses are taken to be PV buses with scheduled voltage at 1.0 p.u.

System operates with the following parameters of interest being observed and evaluated, as shown in Table II.

TABLE II
LINE FLOWS, VOLTAGES & MARGIN INDICES

| Quantity | (p.u) | Quantity | (p.u) |
|---|---|---|---|
| $P_{8-9}$ | 0.8719 | $P_{9-8}$ | -0.8261 |
| $Q_{8-9}$ | 0.1913 | $Q_{9-8}$ | -0.0677 |
| $S_{8-9}$ | 0.8506 | $S_{9-8}$ | -0.8288 |

| Load Bus Number | Index Evaluated | Voltage Magnitude | Voltage angle |
|---|---|---|---|
| 5 | 0.1471 | 0.9834 | -4.587 |
| 7 | 0.1169 | 0.9894 | -1.245 |
| 9 | 0.1957 | 0.9468 | -6.6276 |

It is observed that the load bus 9 has an index close to 0.2 under normal conditions. Now let us consider the loss of (the) line 4 - 9. Table III gives the results for quantity of interest that are necessary to carry out our further discussion.

TABLE II
LINE FLOWS, VOLTAGES & MARGIN INDICES WHEN LINE 4-9 IS DOWN

| Quantity | (p.u) | Quantity | (p.u) |
|---|---|---|---|
| $P_{8-9}$ | 1.7112 | $P_{9-8}$ | -1.4967 |
| $Q_{8-9}$ | 1.6202 | $Q_{9-8}$ | -0.5409 |
| $S_{8-9}$ | 2.3565 | $S_{9-8}$ | 1.5914 |

| Load Bus Number | Index Evaluated | Voltage Magnitude | Voltage angle |
|---|---|---|---|
| 5 | 0.1192 | 0.9752 | -3.2530 |
| 7 | 0.1947 | 0.9283 | -2.3801 |
| 9 | 1.0000 | 0.6147 | -24.6472 |

The results show that under this contingency situation the load bus 9 is very close to steady state voltage collapse situation. The voltage has dropped to a low 0.6147 p.u and the reactive flows and hence line flow has increased significantly which have traditionally been used for sensitivity based voltage collapse detection.

The above simulation brings out the fact that for a security based operating situation for the above test system the normal index at bus 9 has to be kept less than 0.2.

2) Now let us apply the FACTS devices included Voltage Stability constrained OPF, to the WSCC 9 bus test system, to study its effect on load curtailment reduction. For the simulations we have used the following voltage magnitude constraints.

$0.9 \le |V_i| \le 1.1$ for $i$ =1, 2, 3, 4, 6, 8
$0.8 \le |V_i| \le 1.1$ for $i$ =5, 7, 9

Load bus 5 was supposedly having a load demand of 90+ j30 MVA, bus 7 a load demand of 100 + j 35 MVA and load bus 9 having demand of 125 + j 50 MV A. All the generator buses are taken to be PV buses with scheduled voltage at 1.0 p.u. To demonstrate the effectiveness of FACTS devices, let us constrain the load bus 5 with a very strict voltage stability margin of $L_{crit} = 0.1$ under normal conditions. The $L_{crit}$ for load bus 7 and 9 is taken as 0.3. The result of running the voltage stability constrained OPF for various situations of FACTS devices placement and number of FACTS devices is given in Table IV.

TABLE IV
FACTS DEVICES POSITION, CURTAILMENT

| Placement Position | Curtailment at Bus 5 with STATCOM | Curtailment at Bus 5 with SSSC | Curtailment at Bus 5 with UPFC |
|---|---|---|---|
| No FACTS | 0.3592 | 0.3592 | 0.3592 |
| 8-9 | (at 9) 0.4303 | 0.3100 | 0.1867 |
| 5-6 | (at 5) 0.1985 | 0.2081 | 0.0982 |
| 4-5 | (at 4) 0.3635 | 0.2382 | 0.1870 |
| 5-6 & 4-5 | (at 5 & 4) 0.1035 | 0.0890 | 0.0231 |

The strict voltage stability margin index of 0.1 causes load curtailment to have the most effect on load bus 5. However, the incorporation of UPFC causes lowering of the curtailment value. By incorporating STATCOM on bus 5 curtailment reduction is 0.1985 and by SSSC curtailment reduction being 0.281 however as seen from Table IV by incorporating STATCOM in this case is better than SSSC. Placing the FACTS devices on the longest line near the load bus 5 i.e. line 5-6 causes the most efficient reduction when only one TCSC is used. Incorporating two FACTS Devices in lines near bus 5 causes more reduction as the load becomes strongly supported by generators at Bus 1 and Bus 3.

3) Another simulation was carried out to see the impacts of UPFC during contingency situation for different value of stability margin index. The details of the load are the same as presented in sub-section (2). However, the $L_{crit}$ for the load buses 5 and 7 were taken to be 0.3. A contingency of line outage 4-9 was considered. This directly affects the load at bus 9. The results of the simulation without and with UPFC are given in Table V & VI. It was observed that there is no curtailment when UPFC is placed in line 8-9 for this case.

TABLE V
RESULTS WITHOUT UPFC FOR LINE 4-9 OUTAGE

| $L_{crit}$ | Load curtailment at Bus 9 (p.u) | Voltage magnitude at Bus 9 (p.u) |
|---|---|---|
| 0.3 | 0.2842 | 0.9017 |
| 0.4 | 0.1547 | 0.8514 |
| 0.5 | 0.0211 | 0.8041 |

TABLE VI
RESULTS WITH UPFC IN LINE 8-9 FOR LINE 4-9 OUTAGE

| $L_{crit}$ at Bus 9 | Voltage magnitude Bus at 9 (p.u) | Index at Bus 9 |
|---|---|---|
| 0.3 | 0.9871 | 0.3000 |
| 0.4 | 0.9822 | 0.3996 |
| 0.5 | 0.9754 | 0.4387 |

The curtailment was affected when the simulation was carried out without UPFC as shown in Table V. In each of the simulation the $L_{crit}$ constraint caused the load curtailment However, as seen from Table VI by incorporating UPFC, the load curtailment has been avoided even for the strict case of a margin of 0.3. This brings out the fact that UPFC helps in improving the loadability of the system with regards to voltage stability. From a security viewpoint, FACTS Deices helps in maintaining the safety margin for voltage stability margin without compromising on the load.

V. CONCLUSION

The paper discusses the approach to apply the constraint that can take care of incorporating FACTS Deices into the voltage stability constrained OPF algorithm. It is seen that FACTS Deices control improves line flow distribution. It is thus able to reduce load curtailment if any. We have seen from simulations that FACTS Deices relaxes the OPF algorithm when it is constrained by the voltage stability margin indicator. Thus by effectively redistributing the reactive flows, the FACTS Deices aids in relieving voltage stability constrained system operation to an extent.
The applicability of the algorithm in evaluating system reliability measures in composite security based system reliability studies is discussed.

This paper is thus able to formulate an effective method to incorporate FACTS Deices control into an OPF formulation that includes voltage stability. Numerical example illustrates the efficacy of the procedure. The method promises to be a useful tool in security based reliability evaluations of power system operations.

**Ahad Kazemi** He received the M.S degree in electrical engineering from Oklahoma State University, U.S.A., in 1979. Currently, He is an Associated Professor at Electrical Engineering Department of Iran University of Science and Technology, Tehran, Iran. He is the chairman of department power system Iran University of Science and Technology, Tehran, Iran. His research interests are in the application of FACTS Devices in power system control, reactive power planning, and power system restructuring.

**Heidar Ali Shayanfar** He received the B.S. and M.S.E degrees in electrical engineering in 1973 and 1979, respectively. He received the PH. D degree in electrical engineering from Michigan State University, U.S.A., in 1981. Currently, He is a Full Professor at Electrical Engineering Department of Iran University of Science and Technology, Tehran, Iran. His research interests are in the application of artificial intelligence to power system control design, dynamic load modeling, power system absorbability studies and voltage collapse and power system restructuring.

**Masoud Barati** He received the B.S. degree in electrical engineering from Tehran University, Iran. He received the M.S degree in power system from Electrical Engineering Department of Iran University of Science and Technology. His research interests are in the application of FACTS Devices to restructuring power system, Congestion management, power system Dynamics, distributed generation, and power quality.